\documentclass[12pt]{amsart}
\usepackage{amsmath}
\usepackage{amscd}
\usepackage{graphics}
\usepackage{latexsym}

\theoremstyle{plain}

\newtheorem{Thm}{Theorem}
\newtheorem{Lem}[Thm]{Lemma}

\newtheorem{Prop}[Thm]{Proposition}

\theoremstyle{definition}
\newtheorem{Def}{Definition}

\theoremstyle{remark}
\newtheorem{Rem}{Remark}

\def\IC{\mathbb C}

\def\IR{\mathbb C}
\def\IR{\mathbb R}

\def\s2x{\hbox{$S^2 \times S^2$}}
\def\g{\gamma}

\def\bdy{\partial}

	\def\sqr#1#2{{\vcenter{\hrule height.#2pt
    		\hbox{\vrule width.#2pt height#1pt \kern#1pt
       		\vrule width.#2pt}\hrule height.#2pt}}}
	\def\square{\mathchoice\sqr67\sqr67\sqr{2.1}6\sqr{1.5}6}

\begin{document}

\title[]{Lefschetz fibrations on compact Stein surfaces}
\author{Selman Akbulut and Burak Ozbagci} 
\thanks{First named author is partially supported by NSF grant DMS 9971440}
\address{Department of Mathematics \\ Michigan State University \\ MI, 48824}
\email{akbulut@math.msu.edu \and bozbagci@math.msu.edu}
\date{\today}

\begin{abstract}

Let $M$ be a compact Stein surface with boundary. We show that $M$ admits
infinitely many pairwise
nonequivalent positive allowable Lefschetz fibrations over
$D^2$ with bounded fibers.

\end{abstract}
\maketitle
\setcounter{section}{-1}
\section{Introduction}

The existence of a positive allowable
Lefschetz fibration on a
compact Stein surface with boundary
was established by Loi and Piergallini \cite{lp} using
branched covering techniques.
We give an alternative simple  proof
of this fact and construct
explicitly the vanishing
cycles of the Lefschetz fibration, obtaining a direct identification
of compact Stein surfaces with positive allowable Lefschetz fibrations over
$D^2$. In the process we associate to
every compact Stein surface infinitely many pairwise
nonequivalent such Lefschetz
fibrations.

We would like to thank Lee Rudolph, Yasha Eliashberg, Emmanuel Giroux and Ko
Honda for useful discussion about the contact geometry literature.

\section{Preliminaries}

\subsection{Mapping class groups}

Let $F$ be a compact, oriented and connected surface with boundary.
Let $\;\mbox{Diff}^{+} (F, \bdy F)$ be the group of all orientation
preserving self diffeomorphisms of $F$,
fixing  boundary pointwise.
Let $\;\mbox{Diff}_{0}^{+} (F, \bdy F)$ be the subgroup of $\;{\mbox{Diff}}^{+}
(F,\bdy F)$     consisting of all
self diffeomorphisms isotopic to the identity.
Then we define the {\em mapping class group} of the surface $F$ as
$$ Map(F, \bdy F)  = \;{\mbox{Diff}}^{+} (F, \bdy F) / \;{\mbox 
{Diff}}_{0}^{+} (F, \bdy F ).$$

A positive (or right-handed) {\em Dehn twist}
$D(\alpha):F \to F$ about a simple closed curve $\alpha$
is a diffeomorphism obtained by cutting $F$ along $\alpha$,
twisting $360^{\circ}$ to the right and regluing. Note that the positive
Dehn twist $D(\alpha)$ is determined up to isotopy by $\alpha$
and is independent of the orientation on $\alpha$.

It is well-known that the mapping class group
$ Map(F, \bdy F)  $ is generated by Dehn twists. We
will use the functional notation for the products in
$ Map(F, \bdy F) $, e.g., $D(\beta)D(\alpha)$ will
denote the composition where
we apply $D(\alpha)$ first and then $D(\beta)$.

\subsection{ \label{monodromy} Surface bundles over circle }

\indent

In this paper we use the following convention for the monodromy
of a surface bundle
over a circle. We say that an $F$-bundle $W$ over $S^1$
has monodromy $h$ iff $W$ is diffeomorphic to
$$ (F \times I ) / \;(h(x),0) \sim (x,1)$$
where $h \in Map(F, \bdy F) $.
In other words, $h$ is the monodromy if we travel around the
base circle in the positive normal direction to the surface $F$.
Consider the closed 3-manifold
$$ W^{\prime} =  W \cup_{\partial} (\bdy F \times D^2)  .$$
We say that $ W^{\prime} $ has an {\em open book decomposition}
with binding $ \bdy F$, page $F$ and monodromy $h$. It is well-known
that every closed 3-manifold admits an open book decomposition.

\subsection{Positive Lefschetz fibrations}

\indent

Let $M$ be a compact, oriented smooth 4-manifold, and
let $B$ be a compact, oriented smooth 2-manifold. A smooth
map $f:M\to B$ is a {\em positive Lefschetz fibration}
if there exist points $b_1,\ldots,b_m \in $ {\em interior}
$(B)$ such that
\begin{itemize}
\item[(1)] $\{b_1,\ldots,b_m \}$ are the critical values of
$f$, with $p_i\in f^{-1}(b_i)$ a unique critical point
of $f$, for each $i$, and
\item[(2)] about each $b_i$ and $p_i$, there are
local complex coordinate charts agreeing with the orientations
of $M$ and $B$ such that locally $f$ can be expressed
as $f(z_1,z_2)=z_1^2+z_2^2$.
\end{itemize}

It is a consequence of this definition that
$$f|_{f^{-1}(B-\{b_1,\ldots,b_m \})}:
f^{-1}(B-\{b_1,\ldots,b_m \})\to B-\{b_1,\ldots,b_m \}$$
is a smooth fiber bundle over $B-\{b_1,\ldots,b_m \}$
with fiber diffeomorphic to an oriented surface $F$.

Two positive Lefschetz fibrations $f:M\to B$ and
$f^{\prime}:M^{\prime}\to B^{\prime}$ are {\em equivalent}
if there are diffeomorphisms $\Phi:M\to M^{\prime}$ and
$\phi:B\to B^{\prime}$ such that
$f^{\prime}\Phi=\phi f.$

\begin{figure}[ht]
  \begin{center}
     \includegraphics{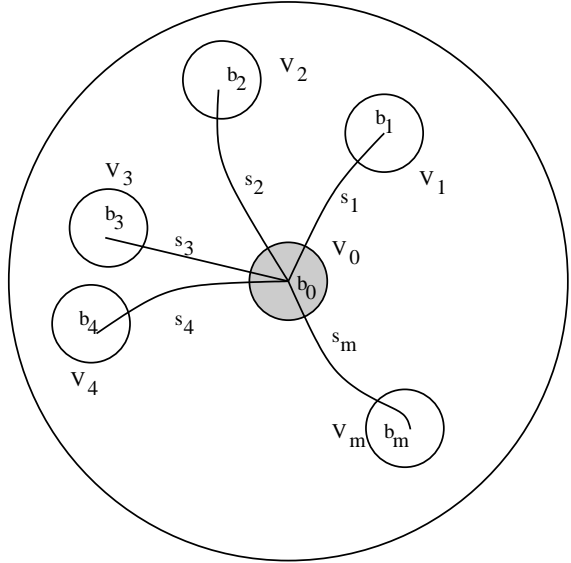}
   \caption{Fibration over the disk} \label{disk}
    \end{center}
  \end{figure}

If $f:M\to D^2$ is a positive Lefschetz fibration,
then we can use this fibration to produce a handlebody
description of $M$.
We select a regular value $b_0 \in$ {\em interior} $(D^2)$
of $f$, an identification $f^{-1}(b_0)\cong F$,
and a collection of arcs $s_i$ in
{\em interior} $(D^2)$ with each $s_i$
connecting $b_0$ to $b_i$, and otherwise disjoint
from the other arcs. We also assume that the critical
values are indexed so that the arcs $s_1,\ldots,s_m$
appear in order as we travel counterclockwise in a
small circle about $b_0$.
Let $V_0,\ldots,V_m$ denote a collection of small disjoint
open disks with $b_i \in V_i$ for each $i$. (cf. Figure~\ref{disk}).

To build our description of $M$, we observe first
that $f^{-1}(V_0)\cong F\times D^2$,
with $\bdy (f^{-1}(V_0)) \cong F\times S^1.$
Let $\nu(s_i)$ be a regular neighborhood of the arc $s_i$.
Enlarging $V_0$ to include the critical value
$b_1$, it can be shown that
$f^{-1}(V_0\cup \nu(s_1) \cup V_1)$
is diffeomorphic to $F\times D^2$
with a 2-handle $H_1$ attached along a circle
${\g_1}$ contained in a fiber
$F\times pt \subset F\times S^1.$
Moreover, condition (2) in the definition of a
Lefschetz fibration requires that $H_1$
is attached with a framing $-1$ relative
to the natural framing on ${\g_1}$
inherited from the product structure of $\bdy (f^{-1}(V_0))$.
${\g_1}$ is called a vanishing cycle.
In addition, $\bdy((F\times D^2)\cup H_1)$ is
diffeomorphic to an $F$-bundle over $S^1$ whose monodromy
is given by $D({\g_1})$, a positive Dehn twist about
${\g_1}$. Continuing counterclockwise about $b_0$,
we add the remaining critical values to our description,
yielding that
$$M_0\cong f^{-1}(
V_0\cup (\bigcup_{i=1}^m \nu(s_i))
\cup (\bigcup_{i=1}^m V_i))$$
is diffeomorphic to $(F\times D^2)\cup
(\bigcup_{i=1}^m H_i)$, where each $H_i$
is a 2-handle attached along a vanishing cycle ${\g_i}$
in an $F$-fiber in $F\times S^1$ with relative framing $-1$.
(For a proof of these statements see \cite{k} or \cite{gs}.)

Furthermore,
$$\bdy M_0 \cong
\bdy((F\times D^2)\cup (\bigcup_{i=1}^m H_i))$$
is an $F$-bundle over $S^1$ with monodromy given
by the composition $D({\g_m})$ $\cdots$ $D({\g_1})$.
We will refer to this product
$D({\g_m})\cdots D({\g_1})$
as the {\em global monodromy} of this fibration.

We note that we can reverse this argument to construct a
positive Lefschetz fibration over $D^2$
from a given set of vanishing cycles.

We say that a positive Lefschetz fibration is {\em allowable} iff all its
vanishing cycles are homologically non-trivial in the fiber $F$.
Note that a simple closed curve on a surface
is homologically trivial iff it separates the surface.

{\Def  PALF is a positive allowable Lefschetz fibration over
$D^2$ with bounded fibers.}

{\Rem \label{palf} With this new notation, we can summarize
the handle attaching procedure
as $$ PALF \; \cup \mbox{ Lefschetz 2-handle} = PALF $$
where a Lefschetz 2-handle is  a
2-handle attached along a nonseparating
simple closed curve in the boundary  with
framing $-1$ relative to the product framing. }

\subsection{ Contact structures}

\indent

We use the standard tight contact structures on ${\IR}^3$, $S^3$  and
${\#}_n S^1 \times S^2 $ (for $n\geq 1$ ) compatible
with their standard orientations.
The structures on  $S^3$  and
${\#}_n S^1 \times S^2 $ are uniquely (up to blowups)
holomorphically fillable ---  $S^3$  as the boundary of
$D^4 \subset {\IC}^2$  and ${\#}_n S^1 \times S^2 $ as the boundary of
$D^4$ union $n$ 1-handles. The tight contact structure on
${\IR}^3 \subset S^3$ will be represented by the kernel of
the 1-form  $ dz +x dy$.

A link $L$ in a contact manifold is called {\em Legendrian}
if its tangent vectors all lie in the contact planes. Legendrian
link theory in ${\IR}^3$ or $ S^3$ reduces to the theory
of the corresponding front projections in ${\IR}^2 $.
We will use projections onto the $yz$-plane in this paper.
The Thurston-Bennequin invariant of a Legendrian knot $L$,
denoted by $tb(L)$, can be computed from a front projection
diagram of $L$ as $$bb(L) - \#  \mbox{left cusps}$$
where $bb(L)$ is the blackboard framing of $L$.

\begin{figure}[ht]
  \begin{center}
     \includegraphics{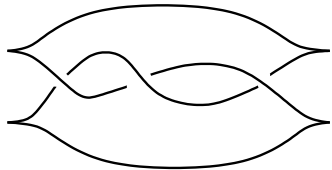}
   \caption{Legendrian trefoil knot} \label{trefoil}
    \end{center}
  \end{figure}

\indent

\section{ Torus knots }

Let $p$ and $q$ be relatively prime integers such that $p,q \geq 2$.

{\Thm \label{torus}

The monodromy of a $(p,q)$ torus knot is a product of
  $(p-1)(q-1)$ nonseparating positive Dehn twists.}

\begin{proof}

It is well-known that a torus knot is  fibered with fiber
being its minimal Seifert surface.
We will describe
how to construct this fiber by plumbing left-handed
Hopf bands (cf. \cite{ha}).

\begin{figure}[ht]
  \begin{center}
     \includegraphics{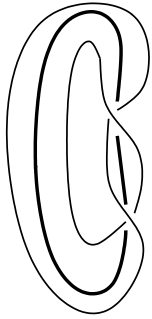}
   \caption{Left-handed Hopf band} \label{leftHopf}
    \end{center}
  \end{figure}

The monodromy of a left-handed Hopf band is
a positive Dehn twist along its core circle as shown in
Figure~\ref{leftHopf}.
Note that our convention for monodromy (see section~\ref{monodromy})
differs from Harer's in  \cite{ha}.

It is proven
in \cite{st} (see also \cite{ga})
that the monodromy of a surface obtained by
plumbing two surfaces is the composition of their
monodromies. We can plumb two
left-handed Hopf bands to get a $(2,3)$ torus knot with
its fibered surface. Simply identify a neigborhood of
the arc $\alpha$ in one
Hopf band with a neighborhood of
the arc $\beta$ in the other Hopf band, transversally as shown in
Figure~\ref{plumbing}. The resulting monodromy will be the product of
two positive Dehn twists along the curves also drawn
in Figure~\ref{plumbing}. Note that the two curves
(one of which is drawn thicker)
intersect each other only once and
they stay parallel when they go
through the left twist on the surface.
It is clear that we can iterate this
plumbing operation to express the monodromy of
a $(2,q)$ torus knot as a product of $(q-1)$ positive Dehn twists.

\begin{figure}[ht]
  \begin{center}
     \includegraphics{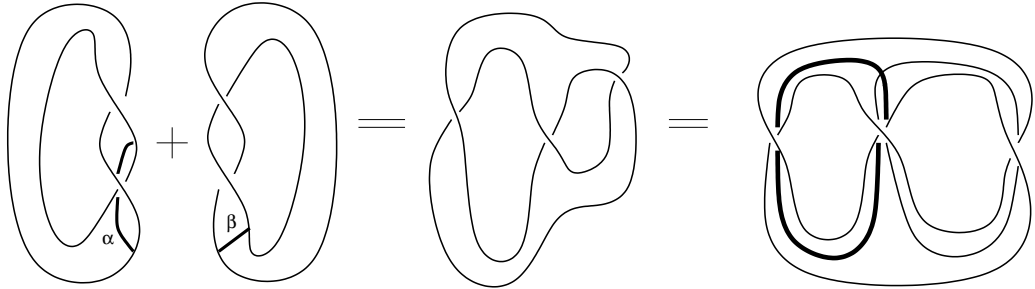}
   \caption{Plumbing two Hopf bands} \label{plumbing}
    \end{center}
  \end{figure}

By attaching more left-handed Hopf bands we can construct
the fibered surface of a $(p,q)$ torus knot for arbitrary $p$ and $q$.
First construct the gate in the back and then plumb a Hopf
band in the front face of that gate and proceed as above to
obtain a second gate. We can iterate this process
to get as many gates as we want.
This is illustrated for $p=3$ and $q=5$ in Figure~\ref{torusknots}.

\begin{figure}[ht]
  \begin{center}
     \includegraphics{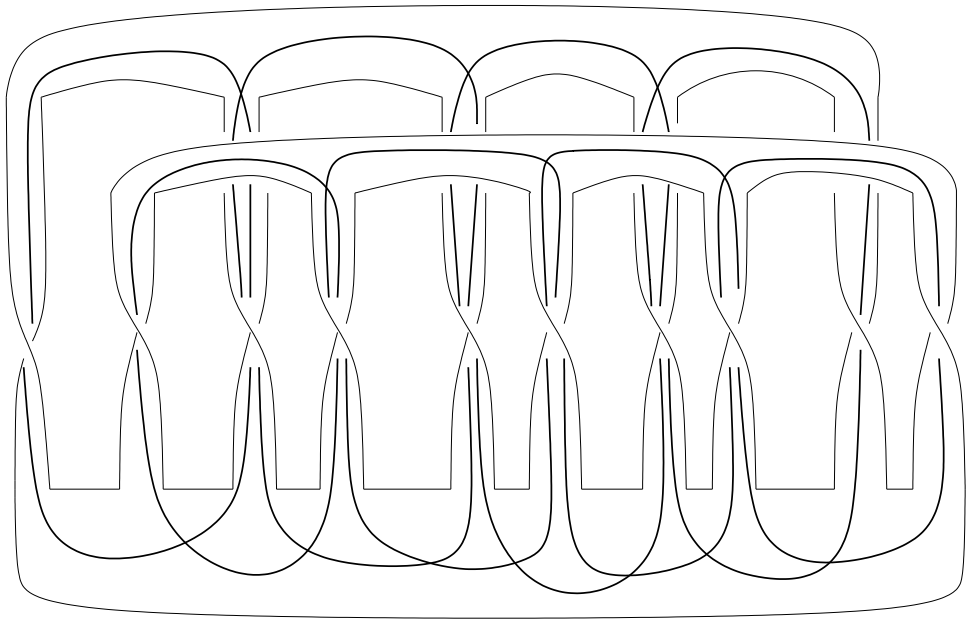}
   \caption{Monodromy of a $(3,5)$ torus knot} \label{torusknots}
    \end{center}
  \end{figure}

Hence the monodromy of the fibration of the complement of a
torus knot in $S^3$ is a product of positive Dehn twists.
These twists are nonseparating by our construction.

\end{proof}

{\Rem Our construction also shows
that the monodromy of an arbitrary
torus link is  a product of positive
Dehn twists. }

{\Thm \label{lyon}  \cite{ly} Let $L$ be a link in $S^3$.
There exists a torus knot $K \subset S^3$ such that
$K \cap L = \emptyset  $ and $L \subset F$
where $F$ is a minimal Seifert surface for $K$. Moreover no component of
$L$ separate the surface $F$.}

\begin{proof}

We describe Lyon's construction given in \cite{ly}.  We say that a link
in $\IR^3$ is in a {\em square bridge position} with respect
to the plane  $ x=0$ if the projection onto the plane
is regular and each segment above the plane projects
to a horizontal segment and each one
below to a vertical segment. Clearly any link
can be put in  a square bridge position.

\begin{figure}[ht]
  \begin{center}
     \includegraphics{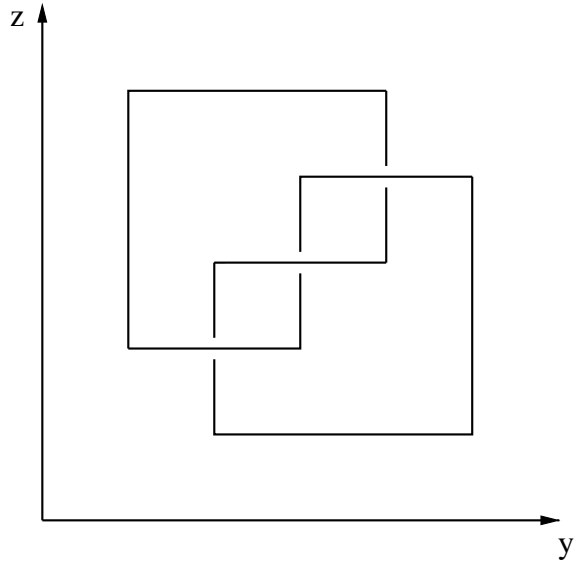}
   \caption{Trefoil knot in a square bridge position} \label{square}
    \end{center}
  \end{figure}

Suppose that
the horizontal and vertical
segments of the projection of the link in the $yz$-plane
are arranged by isotopy so that each horizontal segment is a
subset of $$ \{0 \} \times [0,1] \times \{z_i\} $$
for some $0 < z_1 < z_2 < ... < z_p  < 1 $ and
and each
  vertical
segment is a subset of
$$  \{0 \} \times \{ y_j \}   \times [0,1] $$
for some $ 0 < y_1 < y_2 < ... <y_q < 1$.
Now consider
the 2-disk $$D_i = [\epsilon ,1]  \times [0,1] \times \{z_i\} $$
for each $i=1,2,...,p$
and the 2-disk $$E_j = [ -1, - \epsilon] \times  \{ y_j \}
\times [0,1]$$ for each $j= 1,2,...,q$,
where $\epsilon$ is a small positive number.
Attach these disks by small bands (see Figure~\ref{attaching})
corresponding to each point $(0, y_i , z_j)$
for $i=1,...,p$ and $j=1,...q$. If $p$ and $q$ are
relatively prime then the result is the
minimal Seifert surface $F$ for a $(p,q)$ torus knot $K$
such that $K \cap L = \emptyset  $ and $L \subset F$. Each component of
the link $L$ is a nonseparating embedded curve
on the surface $F$ since we can find an arc
connecting that component to the boundary $K$ from either
side of the component.
Moreover we can choose
$p$ and $q$ arbitrarily large by adding more
disks of either type $D$ or type $E$.

\begin{figure}[ht]
  \begin{center}
     \includegraphics{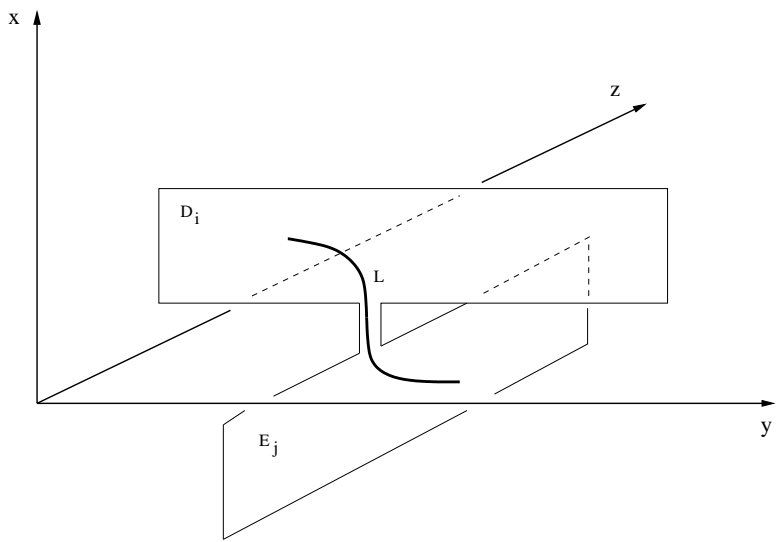}
   \caption{Attaching disks} \label{attaching}
    \end{center}
\end{figure}

\end{proof}

\section{Main theorems}

Let $K$ be a torus knot  in $S^3$. Since $K$ is a fibered knot,
this gives an open book decomposition of $S^3$ with monodromy $h$
which is a product $D({\g_m})\cdots D({\g_1})$
of nonseparating positive Dehn twists given by Theorem~\ref{torus}.
Then $S^3$ bounds a $(PALF)_K$ with global monodromy
$D({\g_m})\cdots D({\g_1})$ and fiber $F$ which is the minimal
Seifert surface for $K$.

{\Prop \label{tor} For any torus knot $K$, $(PALF)_K$  is diffeomorphic to
$D^4$ and has a canonical Stein structure.}

\begin{proof}
Consider the handle decomposition of the $(PALF)_K$ for a
torus knot $K$. Theorem~\ref{torus} gives an explicit
  description of the
vanishing cycles. Cancel each
1-handle with a 2-handle so that the result is just
the 0-handle $D^4$.

\end{proof}

{\Thm (Eliashberg \cite{e}, see also Gompf \cite{g}) A smooth oriented
compact 4-manifold with boundary is a Stein surface,
up to orientation preserving diffeomorphisms, iff
it has a handle decomposition $M_1 \cup H_1 \cup ...
\cup H_n $, where  $M_1$ consists of 0- and 1-handles
and each $H_i$ is a 2-handle attached to $M_1$ along some
attaching circle $L_i$
with framing $tb(L_i) -1 $. }

\indent

We are now ready to state and prove our main theorem.

{\Thm \label{main} Let $M$ be a compact Stein surface with boundary.
Then $M$ admits
infinitely many pairwise
nonequivalent PALF's. Conversely every PALF has a Stein
structure.}

\begin{proof}

Let $M$ be a compact Stein surface with boundary.
We use Eliashberg's characterization of compact
Stein surfaces.

\indent

{\em Case 1 : no 1-handles and one 2-handle }

\indent

Suppose that the compact Stein surface $M$ with
boundary is obtained by attaching
a 2-handle $H$ to $D^4$ along a Legendrian knot
$L$,  with
framing $tb(L)-1 $. Figure~\ref{trefoil} shows
the front projection of a Legendrian trefoil knot.
First of all, we smooth all the cusps of the
diagram and rotate everything counterclockwise to put
$L$ into a square bridge position as in Figure~\ref{square}.

\begin{figure}[ht]
  \begin{center}
     \includegraphics{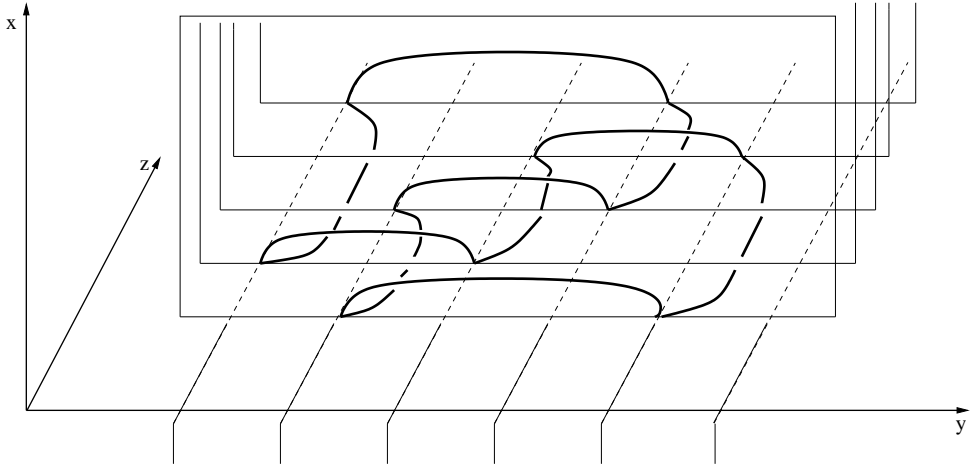}
   \caption{Trefoil knot embedded
into the Seifert surface of a $(5,6)$ torus knot} \label{bridge}
    \end{center}
  \end{figure}

Now we use Lyon's algorithm (cf. Theorem~\ref{lyon})
to find a torus knot
$K$ with its minimal Seifert surface $F$ such that
$L$ is an embedded circle on the
surface $F$. For example, we can embed the
trefoil knot into the
Seifert surface of a (5,6) torus knot as shown in Figure~\ref{bridge}.
Let $L^+$ be a copy of  $L$ pushed in the
positive normal direction to $F$, and let $lk( L, L^+) $ be
the linking number of $L$ and $ L^+ $ computed with
parallel orientations. We need the following observation
to prove our theorem.

{\Lem \label{framing}

$tb(L)  = lk( L, L^+)$ .}

\begin{proof}

When we push $L$ in the positive normal direction to $F$,
we observe that $lk( L, L^+) $ will be exactly the
Thurston-Bennequin framing of $L$, by simply counting the
linking number of $L$ and $L^+$.

\end{proof}

Therefore attaching a 2-handle to $D^4$ along a given
Legendrian knot $L$
in $S^3$, with framing $tb(L) -1 $, is the
same as attaching a 2-handle
along the same knot $L$ (which is isotoped to
be embedded in a fiber of the
boundary of a $(PALF)_K$ ) with framing  $lk( L, L^+) -1 $.
But then the framing $lk( L, L^+) -1 $  is the framing
$-1$ relative to the product framing of $L$.
In other words, we proved that attaching a Legendrian 2-handle is
the same as attaching a Lefschetz 2-handle in our setting.

The global monodromy of $D^4 \cup H \cong (PALF)_K  \cup H$ will
be the monodromy of the torus knot $K$ composed with a positive
Dehn twist along $L$.

\indent

{\em Case 2 :  no 1-handles  }

\indent

Let $L$ be a Legendrian link in $S^3$ with components
$L_1, L_2 , ... , L_n $. Suppose that the compact Stein
surface $M$ with
boundary is obtained by attaching a 2-handle $H_i$
to $D^4$ along $L_i$ for each $ i= 1,2,...,n$.
First
smooth all the cusps of the
diagram and rotate everything counterclockwise to put
$L$ into a square bridge position. Then find a torus knot
$K$ with its minimal Seifert surface $F$
such that each $L_i$ is an embedded circle
on $F$ for $i = 1,2,...,n$.
Now for each $i$, attach a 2-handle
$H_i$ simultaneously to $D^4$
along $L_i$ with framing $lk( L_i, L_i^+) - 1$ .
The result is going to be a PALF by Lemma~\ref{framing} and Remark~\ref{palf},
since the link components
are disjointly embedded nonseparating circles in $F$.

So we showed the global monodromy of $D^4 \cup H_1 \cup ... \cup H_n
\cong (PALF)_K  \cup H_1 \cup ... \cup H_n   $ is the
  monodromy of the torus knot $K$ composed with positive
Dehn twists along $L_i$'s. Note that the Dehn twists along $L_i$'s commute
since they are pairwise disjoint on the surface $F$.

\indent

{\em General case: }

\indent

First we represent the 
1-handles with dotted-circles stacked over 
the front projection of the Legendrian tangle. Here we assume that the 
framed link diagram is in standard form (cf. \cite{g}).
Then we modify the handle decomposition by twisting the strands 
going through each 1-handle negatively once. In the new diagram the Legendrian 
framing will be the blackboard framing with one left-twist added 
for each left cusp. This is illustrated in the second 
diagram in Figure~\ref{tangles}.

\begin{figure}[ht]
  \begin{center}
     \includegraphics{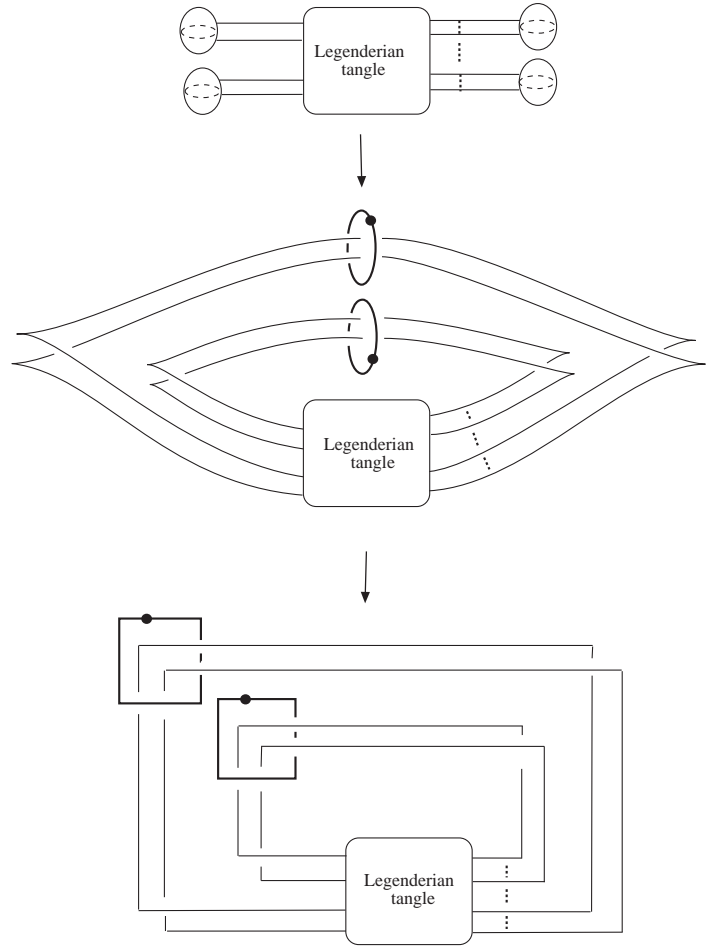}
   \caption{Legendrian link diagram in a 
square bridge position} \label{tangles}
    \end{center}
  \end{figure}

Next we ignore the dots on the dotted-circles for a moment 
and consider the whole diagram as a link in $S^3$. 
Then we put this link diagram in a square bridge position
as in {\em Case 2} (see Figure~\ref{tangles}) 
and find a torus knot $K$ such that
all the link components lie on the
Seifert surface $F$ of $K$. Now consider the $(PALF)_K$ on $D^4$ with 
regular fiber $F$ as in Proposition~\ref{tor}. 
We would like to extend $(PALF)_K$ on $D^4$ to a PALF on $D^4$ 
union 1-handles. 
Recall that attaching a 1-handle to $D^4$ 
(with the dotted-circle notation) is the same as pushing the interior of the 
obvious disk that is spanned by the dotted circle into the interior of 
$D^4$ and removing a tubular neighborhood 
of the image from $D^4$. Before attaching 1-handles we apply 
the following procedure (cf. \cite{ly}): 
We isotope each dotted-circle in the complement of the rest of the 
link such that it becomes transversal to the fibers of
$S^3 \backslash K$, meeting each fiber only once. (see Figure~\ref{bind}).

\begin{figure}[ht]
  \begin{center}
     \includegraphics{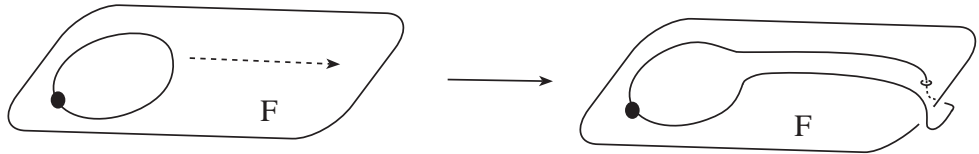}
   \caption{Isotopy of a dotted-circle} \label{bind}
    \end{center}
  \end{figure}

Thus by attaching a 
1-handle to $D^4$ we actually remove a small 2-disk $D^2$ from each fiber 
of $(PALF)_K$, and hence 
obtaining a new PALF on $D^4$ union a 1-handle. 
After attaching all the 1-handles to $D^4$, we get a new PALF 
such that the regular 
fiber is obtained by removing disjoint small disks from $F$. 
Moreover the attaching circles of the 2-handles 
are embedded in a fiber in the boundary of the new PALF such that 
the surface framing of each attaching circle is equal to 
its Legendrian framing. Then {\em Case 2} implies, that is attaching
a Legendrian 2-handle at this stage is the same as attaching
a Lefschetz 2-handle. Hence, 
we can extend our PALF on $D^4 \; \cup$ 1-handles to 
a PALF on $D^4 \; \cup$ 1-handles $\cup$ Legendrian 2-handles. 
The vanishing cycles (hence the monodromy) 
of the constructed PALF are determined explicitly as follows:  
We start with the monodromy of the torus knot $K$, extend 
this over the 1-handles by identity and then we add more 
vanishing cycles 
corresponding to the 2-handles. 

Finally, we note that the $(p,q) $ torus knot in Theorem~\ref{lyon}
can be constructed using arbitrarily large $p$ and $q$.
Therefore our construction
yields infinitely many pairwise
nonequivalent PALF's, since for
chosen $p$ and $q$
the genus of the regular fiber will be at least $(p-1)(q-1)/2$.

Conversely, let $X$ be a PALF, then it is obtained by a
sequence of steps of
attaching 2-handles
$X_{0}= D^2 \times F  \leadsto  X_{1} \leadsto X_{2}.. \leadsto
X_{n}=X \;$, where each $X_{i-1}$ is a PALF
and $X_{i}$ is obtained from  $X_{i-1}$ by attaching a
2-handle to a
nonseparating curve
$C$ lying on a fiber $F\subset \partial X_{i-1}$. Furthermore this handle is
attached to $C$ with the framing $k-1$, where $k$ is the
framing induced from the
surface $F$. Inductively we assume that $ X_{i-1} $ has a Stein
structure, with a
convex fiber $ F \subset \partial X_{i-1}$. By \cite{t} we can start the
induction, and assume that the convex surface $F$ is divided by 
$\partial F$. By
the ``Legendrian realization principle" of
\cite{ho} (pp 323-325), after an isotopy of $(F,C)$,  $k$
can be taken to be the
Thurston-Bennequin framing, and then the result follows by
Eliashberg's theorem (L. Rudolph has pointed out that, in case of $i=1$
identification of $k$ with
Thurston-Bennequin framing also follows from
\cite{r1}- \cite{r4}). Though not necessary, in this process, by 
using \cite{ho}
we can also make the framing of $\partial F$ induced from $F$ to be the
Thurston-Bennequin framing if we wish.

\end{proof}

{\Rem We show in our proof that the PALF structure on a compact 
Stein surface $X$ contains a natural smaller PALF $B^4 \#_{k} S^1 
\times B^3 \to D^2$ given by the associated torus knot, where $k$ 
is the number of 1-handles of $X$.

\vspace{.1in}

{\Rem Our proof shows that by relaxing the condition of positivity,
one can identify smooth bounded $4$-manifolds which are built by $1$- 
and $2$-handles with ALF's (allowable Leschetz fibrations over $D^2$'s). 
In this
case in the proof we start with the binding $ K\; \sharp \;(-K) $ where $K$
is the torus knot, to adjust the framings  (i.e. we use the general form of
\cite{ly}). In particular by \cite{t}, the boundaries of these manifolds
also have contact structures (though not necessarily tight).}

\section{Examples}

\subsection{Example 1}

\begin{figure}[ht]
  \begin{center}
     \includegraphics{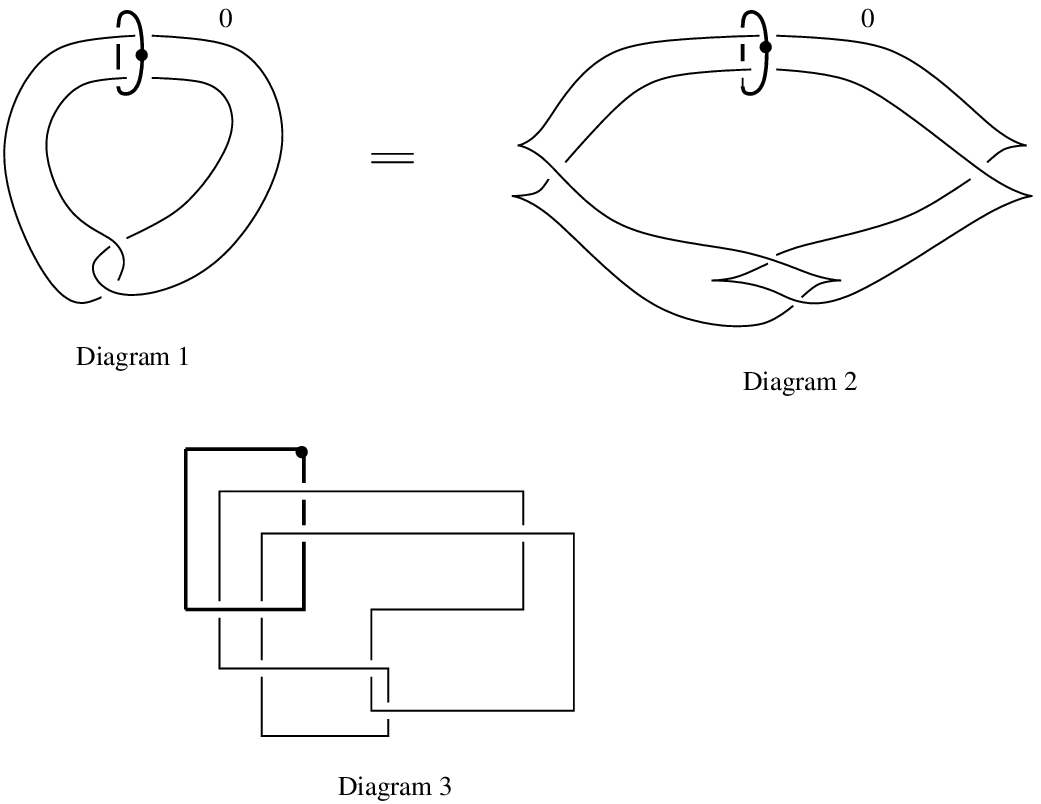}
   \caption{PALF on a fishtail fiber} \label{example1}
    \end{center}
  \end{figure}

In Figure~\ref{example1}, Diagram 1 shows a handle decomposition 
of a smooth 4-manifold $N$ (a regular neighborhood of a fishtail 
fiber in an elliptic fibration) which admits a Stein structure. We modify 
this handle decomposition by twisting the strands going through 
the 1-handle negatively once, as shown in Diagram 2. 
In Diagram 3, we put the whole link (including 
the dotted-circle) in a square bridge position. Note that 
there are exactly $7$ horizontal and $8$ vertical lines in the last diagram. 
Hence according to our algorithm explained above, the Stein surface 
$N$ admits a PALF with $43$ singular fibers where the 
regular fiber is a genus $21$ surface with $2$ boundary components.

\vspace{.1in}

{\Rem The PALF's given by the algorithm of Theorem~\ref{main} 
may not be the most economical ones; 
sometimes with a little care one can find smaller PALF 
structures in the sense of having fewer singular fibers. 
We will illustrate this in the next example.}

\subsection{Example 2}

Let $M$ be the Stein surface given as in Figure~\ref{example2}. 
In the last diagram in Figure~\ref{example2}, 
we put the feet of the 1-handle onto the binding of 
the $(2,2)$ torus link. Since the attaching region (a pair of 3-balls)
of the 1-handle is in a neighborhood of the binding, we can
assume that the pages of the open book will intersect
the pair of balls transversally as 
in Figure~\ref{binding}, so that after gluing the 1-handle to $D^4$
we can extend the fibration over the 1-handle by adding a (2-dimensional)
1-handle to the surface of the fibration without
altering the monodromy. 
Note that this is an alternative way of attaching a 1-handle 
to extend the PALF structure. Hence $M$ admits a PALF
with $2$ singular fibers where the regular fiber is a punctured torus. 
The global monodromy of this PALF is
the monodromy of the $(2,2)$ torus link, extended by identity
over the 1-handle, and composed with a positive Dehn twist
corresponding to the 2-handle. 
In Figure~\ref{fish} we indicate the binding $K$ of the 
open book decomposition of $\bdy M$ obtained from this process.

\begin{figure}[ht]
  \begin{center}
     \includegraphics{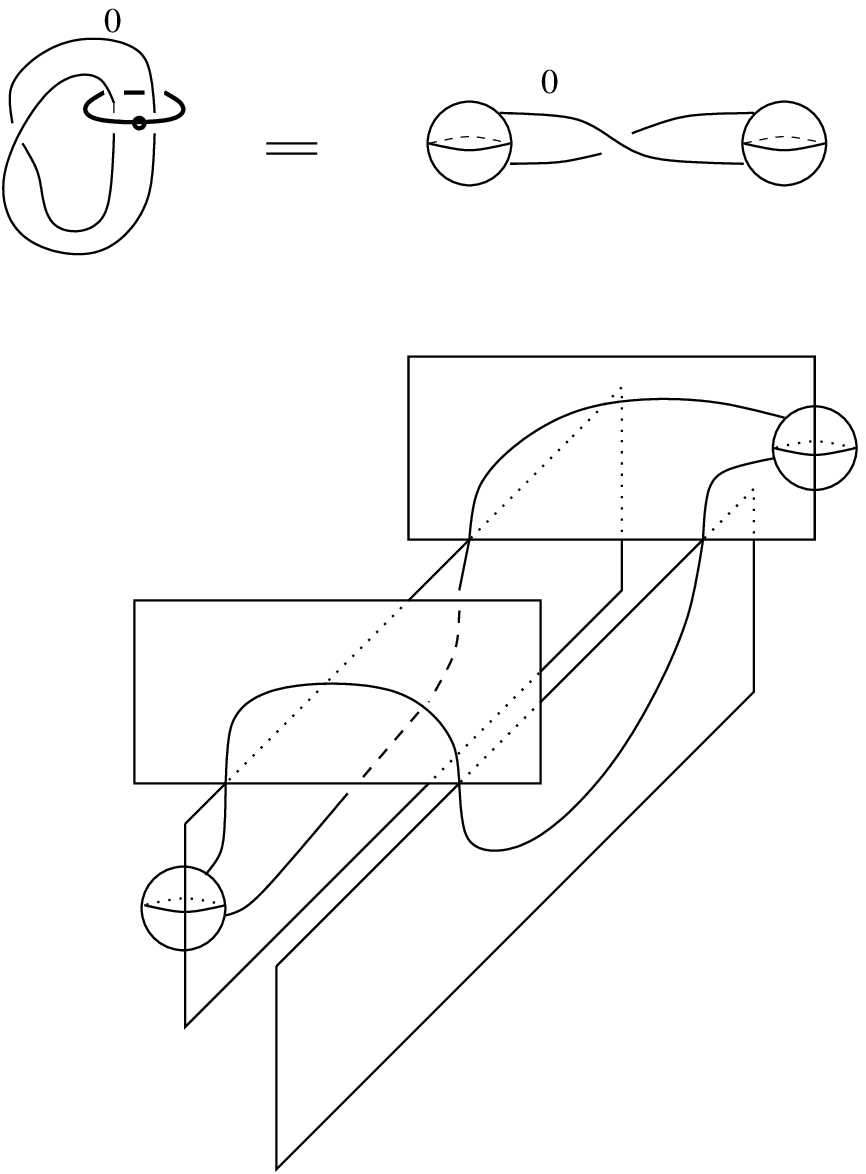}
   \caption{} \label{example2}
    \end{center}
  \end{figure}

\begin{figure}[ht]
  \begin{center}
     \includegraphics{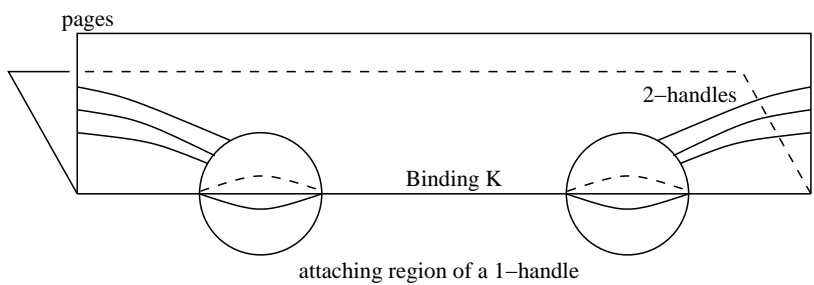}
   \caption{Attaching a 1-handle} \label{binding}
    \end{center}
  \end{figure}

\begin{figure}[ht]
  \begin{center}
     \includegraphics{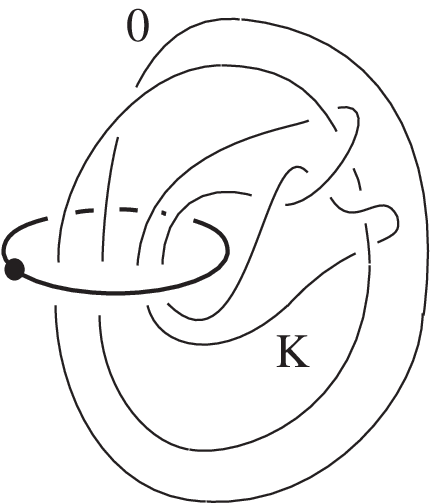}
   \caption{} \label{fish}
    \end{center}
  \end{figure}

\vspace{.1in}

{\Rem In \cite{am} it was shown that every smooth
closed $4$-manifold $X$ can be decomposed as a union of two compact Stein
surfaces along their boundaries $$X=M\cup_{\partial}N. $$
Hence, every $X$ is a union two PALF's along their boundaries. This gives
$4$-manifolds a structure somewhat similar to Heegaard decomposition of $3$
manifolds (we can consider a $3$-dimensional solid handlebody as a
Lefschetz fibration over an interval, with fibers consisting of disks).
Recall that in \cite{am} there is also a relative version of this
theorem; that is, any two smooth closed simply connected
h-cobordant manifolds $ X_{1}, X_{2} $ can be decomposed as union of
Stein surfaces
$ X_{i} = M\cup_{\varphi_{i}}W_{i} $ , where $\varphi_{i}:
\partial W_{i}
\to \partial M $ are diffeomorphisms $ i=1,2 $, $ M $ is simply
connected, and
$W_{1}, W_{2}$ are contractible manifolds which are
diffeomorphic to each other. See also \cite{ao} for more about the 
topology of Stein surfaces. }

\end{document}